\documentstyle[11pt]{article}
\newtheorem{th}{Theorem}
\newtheorem{lem}{Lemma}
\newtheorem{problem}{Problem}
\newcommand{\brm}{\begin{rem}}
\newcommand{\ermq}{\end{rem}}
\newtheorem{pro}{Proposition}
\newtheorem{cor}{Corollary}
\newtheorem{axiom}{Definition}
\newtheorem{rem}{Remark}
\newcommand{\beth}{\begin{th}}
\newcommand{\eeth}{\end{th}}
\newcommand{\bl}{\begin{lem}}
\newcommand{\el}{\end{lem}}
\newcommand{\bp}{\begin{pro}}
\newcommand{\ep}{\end{pro}}
\newcommand{\bcor}{\begin{cor}}
\newcommand{\ecor}{\end{cor}}
\newcommand{\be}{\begin{equation}}
\newcommand{\ee}{\end{equation}}
\newcommand{\beq}{\begin{eqnarray*}}
\newcommand{\eeq}{\end{eqnarray*}}
\newcommand{\beqa}{\begin{eqnarray}}
\newcommand{\eeqa}{\end{eqnarray}}

\newcommand{\cF}{{\cal F}}

\newcommand{\bF}{{\bf F}}

\newcommand{\integ}[2]{\displaystyle \int_{#1}^{#2}}

\newlength{\inter}
\newcommand{\no}{\noindent}
\def \ind{1\!\!1}
\newcommand{\erm}{\end{rem}}
\def \R{I\!\!R}
\def \cadlag {{c\`adl\`ag}~}
\def \esssup {\mbox{ess sup}}
\renewcommand{\baselinestretch}{1.5}
\def \N{I\!\!N}

\title{A Finite Horizon  Optimal Multiple Switching Problem}

\author{Boualem Djehiche,\thanks{Department of Mathematics, The Royal
Institute of Technology, S-100 44 Stockholm, Sweden.\@ e-mail:
boualem@math.kth.se}\,\,\,\, Said Hamad\`ene\thanks{Universit\'e du
Maine, D\'epartement de Math\'ematiques, Equipe Statistique et
Processus, Avenue Olivier Messiaen, 72085 Le Mans, Cedex 9, France.
e-mail: hamadene@univ-lemans.fr}\,\,\, and \, Alexandre
Popier\thanks{Universit\'e du Maine, D\'epartement de
Math\'ematiques, Equipe Statistique et Processus, Avenue Olivier
Messiaen, 72085 Le Mans, Cedex 9, France. e-mail:
Alexandre.Popier@univ-lemans.fr}}
\begin{document}
\date{\today}
\maketitle

%%%%%%%%%%%%%%%%%%%%%Abstract%%%%%%%%%%%%%%%%%%%%%%%%%%

\begin{abstract}
We consider the problem of optimal multiple switching in finite
horizon, when the state of the system, including the switching
costs, is a general adapted stochastic process. The problem is
formulated as an extended impulse control problem and completely
solved using probabilistic tools such as the Snell envelop of
processes and reflected backward stochastic differential equations.
Finally, when the state of the system is a Markov diffusion process,
we show that the vector of value functions of the optimal problem is
a viscosity solution to a system of variational inequalities with
inter-connected obstacles.
\end{abstract}
{\bf AMS Classification subjects}: 60G40 ; 93E20 ; 62P20 ; 91B99.
\medskip

\no {$\bf Keywords$}: real options; backward SDEs; Snell envelope;
stopping time; optimal switching ; impulse control ; variational
inequalities.
\bigskip

%%%%%%%%%%%%%%%%%%%%%%%%%Introduction%%%%%%%%%%%%%%%%%%%%%

\section{\bf Introduction}

Optimal control of multiple switching models arise naturally in many
applied disciplines. The pioneering work by Brennan and Schwartz
(1985), proposing a two-modes switching model for the life cycle of
an investment in the natural resource industry, is probably first to
apply this special case of stochastic impulse control to questions
related to the structural profitability of an investment project or
an industry whose production depends on the fluctuating market price
of a number of underlying commodities or assets. Within this
discipline, Carmona and Ludkosvki (2005) and Deng and Xia (2005)
suggest a multiple switching model to price energy tolling
agreements, where the commodity prices are modeled as continuous
time processes, and the holder of the agreement exercises her
managerial options by controlling the production modes of the
assets. Target tracking in aerospace and electronic systems (cf.\@
Doucet and Ristic (2002)) is another class of problems, where these
models are very useful. These are often formulated as a hybrid state
estimation problem characterized by a continuous time target state
and a discrete time regime (mode) variables. All these applications
seem agree that reformulating these problems in a multiple switching
dynamic setting is a promising (if not the only) approach to fully
capture the interplay between profitability, flexibility and
uncertainty.

\medskip
The optimal two-modes switching problem is probably the most
extensively studied in the literature starting with above mentioned
work by Brennan and Schwartz (1985), and Dixit (1989) who considered
a similar model, but without resource extraction - see Dixit and
Pindyck (1994) and Trigeorgis (1996) for an  overview, extensions of
these models and extensive reference lists. Brekke and \O ksendal
(1991) and (1994), Shirakawa (1997), Knudsen, Meister and Zervos
(1998), Duckworth and Zervos (2000) and (2001), Zervos (2003) and  Pham \& Vath
(2007) use
the framework of generalized impulse control to solve several
versions and extensions of this model, in the case where the
decision to start and stop the production process is done over an
infinite time horizon and the market price process of the underlying
commodity $X$ is a diffusion process, while Trigeorgis (1993) models
the market price process of the commodity as a binomial tree.
Hamad\`ene and Jeanblanc (2007) consider a finite horizon optimal
two-modes switching problem when the price processes are only
adapted to the filtration generated by a Brownian motion while
Hamad\`ene and Hdhiri (2006) extend the set up of the latter paper
to the case where the price processes of the underlying commodities
are adapted to a filtration generated by a Brownian motion and an
independent Poisson process. Porchet {\it et al.} (2006) also study
the same problem, where they assume the payoff function to be given
by an exponential utility function and allow the manager to trade on
the commodities market. Finally, let us mention the work by Djehiche
and Hamad\`ene (2007) where it is shown that including the
possibility of default or bankruptcy in the two-modes switching
model over a finite time horizon, makes the search for an optimal
strategy highly nonlinear and is not at all a trivial extension of
previous results. For example, when the market price of the
underlying commodities is a diffusion process, these optimal
strategies are related to a system of variational inequalities with
inter-connected obstacles, for which very few existence and
regularity results are known in the literature.

An example of the class of multiple switching models discussed in
Carmona and Ludkovski (2005) is related to the management strategies
to run a power plant that converts natural gas into electricity
(through a series of gas turbines) and sells it in the market. The
payoff rate from running the plant is roughly given by the
difference between the market price of electricity and the market
price of gas needed to produce it.

Suppose that besides running the plant at full capacity or keeping
it completely off (the two-modes switching model), there also exists
a total of $q-2$ ($q\geq 3$) intermediate operating modes,
corresponding to different subsets of turbines running.

Let $\ell_{ij}$ denote the switching costs from state $i$ to state
$j$, to cover the required extra fuel and various overhead costs.
Furthermore, let $X=(X_t)_{t\ge 0}$ denote a vector of stochastic
processes that stands for the market price of the underlying
commodities and other financial assets that influence the production
of power. The payoff rate in mode $i$, at time t, is then a function
$\psi_i(t,X_t)$ of $X_t$.

A management strategy for the power plant is a combination of two
sequences:

$(i)$ a nondecreasing sequence of stopping times $(\tau_n)_{n\geq
0}$, where, at time $\tau_n$, the manager decides to switch the
production from its current mode to another one;

$(ii)$ a sequence of indicators $(\xi_n)_{n\geq 1}$ taking values in
$\{1,\ldots,q\}$ of the state the production is switched to. At
$\tau_n$ the station is switched from its current mode $\xi_{n-1}$
to $\xi_n$.
\medskip

When the power plant is run under a strategy $(\delta,
u)=((\tau_n)_{n\geq 1}, (\xi_n)_{n\geq 1})$, over a finite horizon
$[0,T]$, the total expected profit up to $T$ for such a strategy is
$$
J(\delta,u)=E[\integ{0}{T}\psi_{u_s}(s,X_s)ds -\sum_{n\geq
1}\ell_{u_{\tau_{n-1}},u_{\tau_n}}(\tau_n) \ind_{[\tau_{n}<T]}]
$$
where $u_s=\xi_n$ if $s\in [\tau_{n-1},\tau_n[$ $(\tau_0=0)$.

The optimal switching problem we will investigate is to find a
management strategy $(\delta^*, u^*)=((\tau_n^*)_{n\geq
1},(\xi_n^*)_{n\geq 1})$ such that
$$
J(\delta^*, u^*)=\sup_{(\delta,u)} J(\delta,u).
$$

Using purely probabilistic tools such as the Snell envelop of
processes and backward SDEs, inspired by the work by Hamad\`ene and
Jeanblanc (2007), Carmona and Ludkovski (2005) suggest a powerful
robust numerical scheme based on Monte Carlo regressions to solve
this optimal switching problem when $X$ is a diffusion process. They
also list a number of technical challenges, such as the continuity
of the associated value function etc., that prevent a rigorous proof
of the existence and characterization of an optimal solution of this
problem.

The objective of this work is to fill in this gap by providing a
complete treatment of the optimal multiple switching problem, using
the same framework. We are able to prove the existence and provide a
characterization of an optimal strategy of this problem, when $X$
and the switching costs $\ell_{i,j}$ are only adapted to the
filtration $(\cF_t)_{t\ge 0}$ generated by a Brownian motion.

We first provide a Verification Theorem that shapes the problem, via
the Snell envelope of processes. We show that if the Verification
Theorem is satisfied by a vector of continuous processes
$(Y^1,\ldots,Y^q)$ such that, for each $i\in \{1,\ldots,q\}$,
$$
Y^i_t=\esssup_{\tau \geq t}E[\int_t^\tau\psi_i(s,X_s)ds +\max_{j\neq
i}(-\ell_{ij}(\tau)+Y^j_\tau)1_{[\tau <T]}|\cF_t].
$$
then each $Y^i_t$ is the value function of the optimal problem when
the system is in mode $i$ at time $t$:
$$
Y^i_t=\esssup_{(\delta,u)\in {\cal
D}_t}E[\int_t^T\psi_{u_s}(s,X_s)ds- \sum_{j\geq
1}\ell_{u_{\tau_{j-1}} u_{\tau_{j}}}(\tau_j)\ind_{[\tau _j<T]}
|{\cal F}_t].
$$
where ${\cal D}_t$ is the set of strategies such that $\tau_1\geq t$
a.s.\@

An optimal strategy $(\delta^*,u^*)$ is then given by the optimal
stopping times corresponding to the Snell envelop. Moreover, it
holds that $Y^1_0=\sup_{\delta}J(\delta)$, provided that the system
is in mode $i=1$ at time $t=0$.

The unique solution for the Verification Theorem is obtained as the
limit of sequences of processes $(Y^{i,n})_{n\geq 0}$ where,
$Y^{i,n}_t$ is the value function (or the optimal yields) from $t$
to $T$, when the system is in mode $i$ at time $t$ and only at most
$n$ switchings after $t$ are allowed. This sequence of value
functions is defined recursively as follows.
$$
Y^{i,0}_t=E[\int_t^T\psi_i(s,X_s)ds|\cF_t]
$$
and, for $n\ge 1$,
$$
Y^{i,n}_t=\esssup_{\tau \geq t}E[\int_t^\tau\psi_i(s,X_s)ds
+\max_{j\neq i}(-\ell_{ij}(\tau)+Y^{j, n-1}_\tau)1_{[\tau
<T]}|\cF_t].
$$

Finally, if the process $X$ is an It\^o diffusion, with
infinitesimal generator $A$, and each $\ell_{ij}(t)$ is a
deterministic function of $t$, we prove existence of $q$
deterministic continuous functions $v^1(t,x),\ldots,v^q(t,x)$ such
that for any $i\in \{1,\ldots, q\}$, $Y^i_t=v^i(t,X_t)$. Moreover,
the vector $(v^1,\ldots,v^q)$ is a viscosity solution of the
following system of $q$ variational inequalities with
inter-connected obstacles.
$$
\min\{\phi_i(t,x)-\max_{j\neq
i}(-\ell_{ij}(t)+\phi_j(t,x)),-\partial_t\phi_i(t,x)-A\phi_i(t,x)-\psi_i(t,x)\}=0,\,\,\,
\phi_i(T,x)=0, \quad i\in \{1,\ldots,q\}.
$$

\medskip

The organization of the paper is as follows. In Section 2, we give a
formulation of the problem and provide some preliminary results.
Sections 3 \& 4 are devoted to establish the Verification Theorem
and provide an optimal strategy to our problem. In Section
5, we show that, when the driving process $X$ is an It\^o diffusion,
the vector of value functions of our optimal problem is a
viscosity solution of a system of variational inequalities with
inter-connected obstacles. Finally, in Section 6, we provide yet another
numerical scheme
that may be useful in simulating the value-processes satisfying the Verification
Theorem.

%%%%%%%%%%%%%%%%%%%%Section 2 Formulation of the problem
%%%%%%%%%%%%%%%%%%%%%%%%%%%

\section{Formulation of the problem and preliminary results}

The  finite horizon multiple switching  problem can be formulated as
follows. Let ${\cal J}:=\{1,...,q\}$ be the set of all possible
activity modes of the production of the commodity. Being in mode
$i$, a management strategy of the project consists, on the one hand,
of the choice of a sequence of nondecreasing stopping times
$(\tau_n)_{n\geq 1}$ (i.e. $\tau_n\leq \tau_{n+1}$ and $\tau_0=0$)
where the manager decides to switch the activity from its current
mode, $i$, to another one from the set ${\cal
J}^{-i}:=\{1,\ldots,i-1, i+1,\ldots, q\}$. On the other hand, it
consists of the choice of the mode $\xi_n$ to which the production
is switched at $\tau_n$ from the current mode $i$ ; $\xi_n$ is a
random variable which takes its values in ${\cal J}$ and is ${\cal
F}_{\tau_n}$-measurable.

\noindent Assuming that the production activity is in mode $1$ at
the initial time $t=0$, let $(u_t)_{t\leq T}$ denote the indicator
of the production activity's mode at time $t\in [0,T]$:
\begin{equation}
u_t=\ind_{[0,\tau_1]}(t)+\sum_{n\geq1}\xi_n\ind_{(\tau_{n},\tau_{n+1}]}(t).
\end{equation}

\medskip\noindent Note that $\delta:=(\tau_n)_{n\geq 1}$ and the sequence
$\xi:=(\xi_n)_{n\geq 1}$ determine uniquely $u$ and conversely,
$\delta$ and $u$ determine uniquely $(\xi_n)_{n\geq 1}$.

\medskip\noindent A strategy for our multiple switching  problem  will be
simply denoted  by $(\delta,u)$.

\noindent Finally, let $(X_t)_{0\le t\leq T}$  denote the market
price process of e.g.\@ $k$ underlying commodities or other
financial assets that influence the profitability of the production
activity.

\noindent The state of the whole economic system related to the
project at time $t$ is represented by the vector
\begin{equation}
(t, X_t, u_t)\in [0,T]\times \R^k\times {\cal J}.
\end{equation}

\noindent Let $\psi_i(t,x)$ be the payoff rate per unit time when
the system is in state $(t, x, i)$, and for $i,j\in {\cal J}\,\,
(i\neq j)$, $\ell_{ij}:=(\ell_{ij}(t))_{t\leq T}$  denotes the
switching cost of the production at time $t$ from its current mode
$i$ to another mode $j$.

\noindent The expected total profit of running the system with the
strategy $(\delta,u)$ is given by:
$$
J(\delta,u)=E[\integ{0}{T}\psi_{u_s}(s,X_s)ds -\sum_{n\geq
1}\ell_{u_{\tau_{n-1}},u_{\tau_n}}(\tau_n) \ind_{[\tau_{n}<T]}].
$$
Solving the optimal multi-regime starting and stopping problem turns
into finding a strategy $(\delta^*,u^*)$ such that
$J(\delta^*,u^*)\ge J(\delta,u)$ for any other strategies
$(\delta,u)$.

\subsection{Assumptions}\label{assump} Throughout this paper
$(\Omega, {\cal F}, P)$ will be a fixed probability space on which
is defined a standard $d$-dimensional Brownian motion
$B=(B_t)_{0\leq t\leq T}$ whose natural filtration is
$(\cF_t^0:=\sigma \{B_s, s\leq t\})_{0\leq t\leq T}$. Let $
\bF=(\cF_t)_{0\leq t\leq T}$ be the completed filtration of
$(\cF_t^0)_{0\leq t\leq T}$ with the $P$-null sets of ${\cal F}$.
Hence $\bF$ satisfies the usual conditions, $i.e.$, it is right
continuous and complete.

Furthermore, let:
\begin{itemize}
\item[-] ${\cal P}$ be the $\sigma$-algebra on $[0,T]\times \Omega$ of
$\bF$-progressively measurable sets ;
\item[-] ${\cal M}^{p,l}$
be the set of $\cal P$-measurable and $\R^l$-valued processes
$w=(w_t)_{t\leq T}$ such that $E[\int_0^T|w_s|^pds]<\infty$  and
${\cal S}^p$  be the set of $\cal P$-measurable, continuous,
$\R$-valued processes ${w}=({w}_t)_{t\leq T}$ such that
$E[\sup_{0\le t\leq T}|{w}_t|^p]<\infty$ ($p>1$ is fixed) ;
\item[-] For any stopping time $\tau \in [0,T]$, $\,{\cal T}_\tau$ denotes
the set of all stopping times $\theta$ such that $\tau \leq \theta
\leq T$, $P-a.s.$.
\end{itemize}

We now make the following assumptions on the data:
\begin{itemize}
\item[(i)] The market price $X:=(X_t)_{0\leq
t\leq T}$ is $\R^k$-valued and each component belongs to ${\cal
S}^p$.
\item[(ii)] The functions $\psi_i(t,x)$, $(t,x)\in [0,T]\times
\R^k$ and $i\in {\cal J}$, are continuous and satisfy a linear
growth condition, i.e. there exists a constant $C$ such that
$|\psi_i(t,x)|\leq C(1+|x|)$ for $0\leq t\leq T$ and $x\in \R^k$.
\item[(iii)] The processes $\ell_{ij}$ belong to ${\cal S}^p$ and
there exists a real constant $\gamma>0$ such P-a.s. for any $0\leq
t\leq T$, $\min\{\ell_{ij}(t), i,j\in {\cal J}, \,i\neq j\}\geq
\gamma$.
\item[(iv)] $(\tau_n)_{n\geq 1}$ are $\bF$-stopping times and $(\xi_n)_{n\geq
1}$ are random variables with values in ${\cal J}$ and such that for any
$n\geq 1$, $\xi_n$ is ${\cal F}_{\tau_n}$-measurable. Additionally,
we assume that for any $n\geq 1$, $P[\xi_n=\xi_{n+1}]=0$. The
strategies $(\delta , u)=((\tau_n)_{n\geq 1},(\xi_n)_{n\geq 1})$ are
called admissible if they satisfy:
$$
\lim_{n\to \infty}\tau_n =T\quad P-\mbox{a.s.}
$$
The set of admissible strategies is denoted by ${\cal D}_a$.
\end{itemize}
\brm \label{remq1} The above assumptions on $X$ and $\psi_i$,
i=1,\ldots,q, can be modified or weakened in any way which preserves
the fact that the process $(\psi_i(t,X_t); \,\, 0\le t\leq T,\,\,
i\in {\cal J})$ belongs to ${\cal M}^{p,1}$.\erm

We can now formulate the multi-regime starting and stopping problem
as follows:
\begin{problem}\label{optimal} Find a strategy
$(\delta^*,u^*)=((\tau_n^*)_{n\geq 1},(\xi_n^*)_{n\geq 1})\in {\cal
D}_a$ such that
\begin{equation}
J(\delta^*,u^*)=\sup_{(\delta ,u)\in {\cal D}_a} J(\delta,u).
\end{equation}
\end{problem}

\noindent An admissible strategy $(\delta,u)$ is called {\it finite}
if, during the time interval $[0,T]$, it allows the manager to make
only a finite number of decisions, i.e. $P[\omega,
\tau_n(\omega)<T,\,\,\, \mbox{for all}\,\,\, n\geq 0]=0$. Hereafter
the set of finite strategies will be denoted by ${\cal D}$. The next
proposition tells us that the supremum of the expected total profit
can only be reached over finite strategies .

\bp \label{finite} The suprema over admissible strategies and finite
strategies coincide:
\begin{equation}
\sup_{(\delta,u) \in {\cal D}_a}J(\delta,u)=\sup_{(\delta,u) \in
{\cal D}}J(\delta,u).
\end{equation}
\ep {\it Proof}. If  $(\delta,u)$ is an admissible strategy which
does not belong to ${\cal D}$, then $J(\delta)=-\infty$. Indeed, let
$A=\{\omega, \tau_n(\omega)<T, \,\,\, \mbox{for all}\,\,\, n\geq
0\}$ and $A^c$ be its complement. Since $(\delta,u) \in {\cal D}_a
\setminus {\cal D}$, then $P(A)>0$. Since  the process $X$ belongs
to ${\cal S}^{p}$ and $\psi_i$ is of linear growth, then the
processes $(\psi_i(t,X_t))_{t\leq T}$ belongs to ${\cal M}^{p,1}$.
Therefore,
$$
\begin{array}{ll}
J(\delta,u)& \leq E[\integ{0}{T}(\max_{i\in {\cal J}}{
|\psi_i(s,X_s)|})\,\, ds]\\& -E[\{\sum_{n\geq
1}\ell_{u_{\tau_{n-1}},u_{\tau_n}}(\tau_n) \} \ind_A +\{\sum_{n\geq
1}\ell_{u_{\tau_{n-1}},u_{\tau_n}}(\tau_n)
\ind_{[\tau_{n}<T]}]\}\ind_{A^c}]=-\infty,\end{array}
$$
since for any $t\leq T$ and $i,j\in {\cal J}$, $\ell_{ij}(t)\geq
\gamma 0$. This implies that $J(\delta,u)=-\infty$ and then\\
$\sup_{(\delta,u) \in {\cal D}_a}J(\delta,u)=\sup_{(\delta,u) \in
{\cal D}}J(\delta,u)$. $\Box$
\bigskip

We finish this section by introducing the key ingredient of the
proof of the main result, namely the notion of Snell envelope and
its properties. We refer to Cvitanic and Karatzas (1996) , Appendix
D in Karatzas and Shreve (1998), Hamad\`ene (2002) or El Karoui
(1980) for further details.

\subsection{The Snell Envelope}

In the following  proposition we summarize the main results on the
Snell envelope of processes used in this paper.

\bp \label{thmsnell} Let $U=(U_t)_{0\le t\leq T}$ be an
$\bF$-adapted $\R$-valued \cadlag process that belongs to the class
[D], $i.e.$\@ the set of random variables $\{U_\tau, \,\,\tau
\in{\cal T}_0\}$ is uniformly integrable. Then, there exists an
$\bF$-adapted $\R$-valued \cadlag process $Z:=(Z_t)_{0\le t\leq T}$
such that:
\begin{itemize}

\item[(i)] $Z$ is the smallest super-martingale which dominates $U$,
$i.e$, if $(\bar{Z}_t)_{0\leq t\leq T}$ is another \cadlag
supermartingale of class [D] such that for all $0\leq t\leq T$,
$\bar{Z}_t\geq U_t$ then $\bar{Z}_t\geq Z_t$ for any $0\leq t\leq
T$.

\item[(ii)] For any $\bF$-stopping time $\theta$ we have:  \be
\label{sun} Z_\theta=\esssup_{\tau \in {\cal
T}_{\theta}}E[U_\tau|\cF_\theta]\,\,\,\,\,\,(\mbox{and then
}Z_T=U_T).\ee The process $Z$ is called the {\it Snell envelope\,}
of $U$.
\end{itemize}

Moreover, the following properties hold:
\begin{itemize}
\item[(iii)] The Dood-Meyer decomposition of $Z$ implies the
existence of a martingale $(M_t)_{0\le t\leq T}$ and two
nondecreasing processes $(A_t)_{0\le t\leq T}$ and $(B_t)_{0\le
t\leq T}$ which are respectively continuous and purely discontinuous
predictable such that for all $0\le t\leq T$,
$$
Z_t=M_t-A_t-B_t \quad (\mbox{with } A_0=B_0=0).
$$
Moreover, for any $0\le t\leq T$, $\{\Delta_tB >0\}\subset \{\Delta
_tU<0\}\cap \{Z_{t-}=U_{t-}\}.$
\item[(iv)] If $U$ has only positive jumps then $Z$ is a continuous process.
Furthermore, if $\theta$ is an $\bF$-stopping time and
$\tau^*_\theta=\inf\{s\geq \theta, Z_s=U_s\}\wedge T$ then
$\tau^*_\theta$ is optimal after $\theta$, $i.e.$\@
\begin{equation}\label{sdeux}Z_\theta=E[Z_{\tau^*_\theta}|\cF_\theta]=E[U_{\tau^*_\theta}|\cF_\theta]=\esssup
_{\tau \geq \theta}E[U_\tau|\cF_\theta].\end{equation}
\item[(v)] If $(U^n)_{n\geq 0}$ and $U$ are \cadlag and of class [D] and such
that the sequence $(U^n)_{n\geq 0}$ converges increasingly and
pointwisely to $U$ then $(Z^{U^n})_{n\geq 0}$ converges increasingly
and pointwisely to $Z^U$; $Z^{U_n}$ and $Z^U$ are the Snell
envelopes of respectively $U_n$ and $U$. Finally, if $U$ belongs to
${\cal S}^p$ then $Z^U$ belongs to ${\cal S}^p$.
\end{itemize}
\ep

For the sake of completeness, we give a proof of the stability
result $(v)$.

\noindent {\it Proof of}\,\,$(v)$. Since, for any $n\geq 0$, $U^n$
converges increasingly and pointwisely to $U$, it follows that for
all $t\in [0,T]$, $Z^{U_n}_t\leq Z^{U}_t$ $P$-a.s. Therefore,
$P-a.s.$, for any $t\in [0,T]$, $\lim_{n\rightarrow
\infty}Z^{U_n}_t\leq Z^{U}_t$. Note that the process
$(\lim_{n\rightarrow \infty}Z^{U^n}_t)_{0\le t\leq T}$ is a \cadlag
supermartingale of class [D], since it is a limit of a nondecreasing
sequence of supermartingales (see e.g. Dellacherie and Meyer (1980),
pp.86). But $U^n\leq Z^{U_n}$ implies that $P-a.s.$, for all $t\in
[0,T]$, $U_t\leq \lim_{n\rightarrow \infty}Z^{U^n}_t$. Thus,
$Z^{U}_t\leq \lim_{n\rightarrow \infty}Z^{U^n}_t$ since the Snell
envelope of $U$ is the lowest supermartingale that dominates $U$. It
follows that $P-a.s.$, for any $t\leq 1$, $\lim_{n\rightarrow
\infty}Z^{U^n}_t=Z^{U}_t$, whence the desired result.

Assume now that $U$  belongs to ${\cal S}^p$. Since, for any $0\le
t\leq T$, $-E[\sup_{0\le s\leq T}|U_s||{\cal F}_t]\leq U_t\leq
E[\sup_{0\le s\leq T}|U_s||{\cal F}_t]$, using the Doob-Meyer
inequality, it follows that $Z^U$ also belongs to ${\cal S}^p$.
$\Box$

%%%%%%%%%%%%%%%%verification Theorem%%%%%%%%%%%%%%%%%%%%%%%%%%%%%%%%%

\section{ A verification Theorem}

In terms of a verification theorem,  we show that Problem
\ref{optimal} is reduced to the existence of  $q$ continuous
processes $Y^1,\ldots,Y^q$ solutions of a system of equations
expressed via the Snell envelopes. The process $Y^i_t$, for $i\in
{\cal J}$, will stand for the optimal expected profit if, at time
$t$, the production activity is in the state $i$.

\noindent Set
$$
D_\tau(\zeta=\zeta'):= \inf\{s\geq \tau, \zeta_s=\zeta'_s\}\wedge T,
$$
where, $\tau$ is an $\bf F$-stopping time and $(\zeta_t)_{0\leq
t\leq T}$, $(\zeta'_t)_{0\le t\leq T}$ are two continuous $\bf
F$-adapted and $\R-$valued processes.

\medskip

\beth \label{thmverif}$($Verification Theorem$)$

\noindent Assume there exist $q$ ${\cal S}^p$-processes
$(Y^i:=(Y^i_t)_{0\le t\leq T},\,\,i=1,\ldots,q)$ that satisfy \be
\begin{array}{l}
\label{eqvt} Y^i_t=\esssup_{\tau \geq t}E[\int_t^\tau\psi_i(s,X_s)ds
+\max_{j\in {\cal J}^{-i}}(-\ell_{ij}(\tau)+Y^j_\tau)1_{[\tau
<T]}|\cF_t]\,\quad (\mbox{and then } Y^i_T=0).
\end{array}
\ee Then $Y^1,\dots,Y^q$ are unique. Furthermore :
\begin{itemize}
\item[$(i)$]
\be\label{Yzero} Y^1_0=\sup_{(\theta,v) \in {\cal D}}J(\theta,v). \ee
\item[$(ii)$] Define the sequence
of $\bF$-stopping times $(\tau_n)_{n\geq 1}$ by \be\label{tau1}
\begin{array}{l}
\tau_1=D_0(Y^1= \max_{j\in {\cal J}^{-1}}(-\ell_{1j}+Y^j))
\end{array}
\ee and, for $n\geq 2$, \be\label{taun}
\begin{array}{l}
\tau_n=D_{\tau_{n-1}}(Y^{u_{\tau_{n-1}}}=\max_{k\in {\cal
J}^{-\tau_{n-1}}}(-\ell_{\tau_{n-1}k}+Y^k)),
\end{array}
\ee where,

$\bullet$ $u_{\tau_1}=\sum_{j\in {\cal J}}j\ind_{\{\max_{k\in {\cal
J}^{-1}}(-\ell_{1k}({\tau_1})+Y^k_{\tau_1})=-\ell_{1j}(\tau_1)+Y^j_{\tau_1}\}};$

$\bullet$ For any $n\geq 1$ and $t\geq \tau_n$,
$Y_t^{u_{\tau_{n}}}=\sum_{j\in {\cal J}}\ind_{[u_{\tau_n}=j]}Y^j_t$;

$\bullet$ For $n\geq 2$, $u_{\tau_{n}}=l$ on the set $\{\max_{k\in
{\cal
J}^{-u_{\tau_{n-1}}}}(-\ell_{u_{\tau_{n-1}}k}({\tau_{n}})+Y^k_{\tau_{n}})
=-\ell_{u_{\tau_{n-1}}l}(\tau_{n})+Y^l_{\tau_{n}}\}$,

where, $\ell_{u_{\tau_{n-1}}k}({\tau_{n}})=\sum_{j\in {\cal
J}}\ind_{[\tau_{n-1}=j]}\ell_{jk}({\tau_{n}})$ and ${\cal
J}^{-u_{\tau_{n-1}}}=\sum_{j\in {\cal J}}\ind_{[\tau_{n-1}=j]}{\cal
J}^{-j}$.
\medskip

Then, the strategy $(\delta,u)$ is optimal i.e.\@ $J(\delta,u)\ge
J(\theta,v)$ for any $(\theta,v) \in {\cal D}$.
\end{itemize}

\eeth \no $Proof$. The proof consists essentially in showing that
each process $Y^i$, as defined by (\ref{eqvt}), is nothing but the
expected total profit or the value function of the optimal problem,
given that the system is in mode $i$ at time $t$. More precisely,
\be\label{eqyn} Y^i_t=\esssup_{(\delta,u)\in {\cal
D}_t}E[\int_t^T\psi_{u_s}(s,X_s)ds- \sum_{j\geq
1}\ell_{u_{\tau_{j-1}} u_{\tau_{j}}}(\tau_j)\ind_{[\tau _j<T]}
|{\cal F}_t], \ee where ${\cal D}_t$ is the set of finite strategies
such that $\tau_1\geq  t$, $P$-a.s.\@ if at time $t$ the system is
in the mode $i$. This characterization implies in particular that
the processes $Y^1,\ldots,Y^q$ are unique. Moreover, thanks to a
repeated use of the characterization of the Snell envelope
(Proposition \ref{thmsnell}, $(iv)$), the strategy $(\delta, u)$
defined recursively by (\ref{tau1}) and (\ref{taun}), is shown to be
optimal.

\noindent Indeed, since at time $t=0$ the system is in mode $1$, it
holds true that, for any $ 0\le t\le T$, \be
\label{eq1}Y^1_t+\integ{0}{t} \psi_1(s,X_s) ds= \esssup _{\tau \geq
t}E[\integ{0}{\tau} \psi_1(s,X_s)ds+\max_{j\in {\cal
J}^{-1}}(-\ell_{1j}(\tau)+Y^j_\tau)\ind_{[\tau <T]}|\cF_t]. \ee But,
$Y^1_0$ is $\cF_0$-measurable. Therefore it is $P-a.s.$ constant and
then $Y^1_0=E[Y^1_0]$.

\noindent On the other hand, according to Proposition
\ref{thmsnell}, $(iv)$, $\tau_1$ as defined by (\ref{tau1}) is
optimal,  $Y^1_T=0$ and
$$
u_{\tau_1}=\sum_{j\in {\cal J}}j\ind_{\{\max_{k\in {\cal
J}^{-1}}(-\ell_{1k}({\tau_1})+Y^k_{\tau_1})=-\ell_{1j}(\tau_1)+Y^j_{\tau_1}\}}.
$$
Therefore, \be \label{eq2}
\begin{array}{ll}
Y^1_0&=E[\integ{0}{\tau_1} \psi_1(s,X_s)ds+\max_{j\in {\cal
J}^{-1}}(-\ell_{1j}(\tau_1)+Y^j_{\tau_1})\ind_{[\tau _1<T]}]\\
{}&=E[\integ{0}{\tau_1} \psi_1(s,X_s)ds+(-\ell_{1
u_{\tau_1}}(\tau_1)+Y^{u_{\tau_1}}_{\tau_1})\ind_{[\tau
_1<T]}].\end{array}\ee

\noindent Next, we claim that  $P$-a.s.\@  for every $t\in
[\tau_1,T],$ \be \label{eq3}Y^{u_{\tau_1}}_{t}=\esssup_{\tau \geq
t}E[\int_{t}^{\tau}\psi_{u_{\tau_1}}(s,X_s)ds +\max_{j\in {\cal
J}^{-{u_{\tau_1}}}}(-\ell_{u_{\tau_1}j}(\tau)+Y^j_\tau)1_{[\tau
<T]}|\cF_{t}]. \ee To see this, recall that for any $i\in{\cal J}$
and  $0\le t\leq T$
$$
Y^i_t=\esssup_{\tau \geq t}E[\int_t^\tau\psi_i(s,X_s)ds +\max_{j\in
{\cal J}^{-i}}(-\ell_{ij}(\tau)+Y^j_\tau)1_{[\tau <T]}|\cF_t].
$$
This also means that the process
$\left(Y^i_t+\int_0^t\psi_i(s,X_s)ds\right)_{0\le t\leq T}$ is a
supermartingale which dominates
$$
\left(\int_0^t\psi_i(s,X_s)ds +\max_{j\in {\cal
J}^{-i}}(-\ell_{ij}(t)+Y^j_t)1_{[t <T]}\right)_{0\le t\leq T}.
$$
This implies that the process
$\left(\ind_{[u_{\tau_1=i]}}(Y^i_t+\int_{{\tau_1}}^t\psi_i(s,X_s)ds)\right)_{t\in
[\tau_1,T]}$ is a supermartingale which dominates
$$
\left(\ind_{[u_{\tau_1=i]}}(\int_{\tau_1}^t\psi_i(s,X_s)ds
+\max_{j\in {\cal J}^{-i}}(-\ell_{ij}(t)+Y^j_t)1_{[t
<T]})\right)_{t\in [\tau_1,T]}.
$$
Since ${\cal J}$ is finite, the process $\left(\sum_{i\in {\cal
J}}\ind_{[u_{\tau_1=i]}}(Y^i_t+\int_{{\tau_1}}^t\psi_i(s,X_s)ds)\right)_{t\in
[\tau_1,T]}$ is also a supermartingale which dominates
$\left(\sum_{i\in {\cal
J}}\ind_{[u_{\tau_1=i]}}(\int_{\tau_1}^t\psi_i(s,X_s)ds +\max_{j\in
{\cal J}^{-i}}(-\ell_{ij}(t)+Y^j_t)1_{[t <T]})\right)_{t\in
[\tau_1,T]}$.

\noindent Thus, the process
$\left(Y^{u_{\tau_1}}_t+\int_{{\tau_1}}^t\psi_{u_{\tau_1}}(s,X_s)ds\right)_{t\in
[\tau_1,T]}$ is a supermartingale which is greater than \be
\left(\int_{\tau_1}^t\psi_{u_{\tau_1}}(s,X_s)ds +\max_{j\in {\cal
J}^{-{u_{\tau_1}}}}(-\ell_{{u_{\tau_1}}j}(t)+Y^j_t)1_{[t
<T]}\right)_{t\in [\tau_1,T]}. \ee

\noindent To complete the proof it remains to show that it is the
smallest one which has this property and use the characterization of
the Snell envelope (Proposition \ref{thmsnell}, $(i)-(ii)$).

\noindent Indeed, let $(Z_t)_{t\in [\tau_1,T]}$ be a supermartingale
of class [D] such that, for any $t\in [\tau_1,T]$,
$$
Z_t\geq \int_{\tau_1}^t\psi_{u_{\tau_1}}(s,X_s)ds +\max_{j\in {\cal
J}^{-{u_{\tau_1}}}}(-\ell_{{u_{\tau_1}}j}(t)+Y^j_t)1_{[t <T]}.
$$
It follows that for every  $t\in [\tau_1,T]$,
$$
Z_t \ind_{[u_{\tau_1}=i]}\geq
\ind_{[u_{\tau_1}=i]}\left(\int_{\tau_1}^t\psi_{i}(s,X_s)ds +
\max_{j\in{\cal J}^{-{i}}}(-\ell_{ij}(t)+Y^j_t)1_{[t <T]}\right).
$$
But, the process $(Z_t \ind_{[u_{\tau_1}=i]})_{t\in [\tau_1, T]}$ is
a supermartingale and for every  $t\in [\tau_1,T]$,
$$
\ind_{[u_{\tau_1}=i]}Y^i_t=\esssup_{\tau \geq
t}E[\ind_{[u_{\tau_1}=i]}(\int_t^\tau\psi_i(s,X_s)ds +\max_{j\in
{\cal J}^{-i}}(-\ell_{ij}(\tau)+Y^j_\tau)1_{[\tau <T]})|F_t].
$$
%Actually, we can put the indicator
%$\ind_{[u_{\tau_1}=i]}$ inside the $\esssup$ since it is ${\calF}_t$-measurable ($t\geq
%\tau_1$) and non-negative.
It follows that, for every  $t\in [\tau_1,T]$,
$$
\ind_{[u_{\tau_1}=i]}Z_t\geq
\ind_{[u_{\tau_1}=i]}(Y^i_t+\int_{\tau_1}^t\psi_{i}(s,X_s)ds).
$$
Summing over $i$, we get, for every  $t\in [\tau_1,T]$,
$$
Z_t\geq Y^{u_{\tau_1}}_t+\int_{\tau_1}^t\psi_{u_{\tau_1}}(s,X_s)ds.
$$
Hence, the process
$\left(Y^{u_{\tau_1}}_t+\int_{\tau_1}^t\psi_{u_{\tau_1}}(s,X_s)ds\right)_{t\in
[\tau_1,T]}$ is the Snell envelope of
$$
\left(\int_{\tau_1}^t\psi_{u_{\tau_1}}(s,X_s)ds +\max_{j\in {\cal
J}^{-{u_{\tau_1}}}}(-\ell_{{u_{\tau_1}}j}(t)+Y^j_t)1_{[t
<T]}\right)_{t\in [\tau_1,T]},
$$
whence Eq. (\ref{eq3}).

\medskip\noindent
Now, from (\ref{eq3}) and the definition of $\tau_2$ in Eq.\@
(\ref{taun}), we have
$$
\begin{array}{ll}
Y^{u_{\tau_1}}_{\tau_1}&=E[\int_{\tau_1}^{\tau_2}\psi_{u_{\tau_1}}(s,X_s)ds
+\max_{j\in {\cal
J}^{-{u_{\tau_1}}}}(-\ell_{u_{\tau_1}j}(\tau_2)+Y^j_{\tau_2})1_{[\tau_2
<T]}|F_{\tau_1}]\\
{}&=E[\int_{\tau_1}^{\tau_2}\psi_{u_{\tau_1}}(s,X_s)ds
+(-\ell_{u_{\tau_1}{u_{\tau_2}}}(\tau_2)+Y^{u_{\tau_2}}_{\tau_2})1_{[\tau_2
<T]}|F_{\tau_1}].\end{array}
$$
Setting this characterization of $Y^{u_{\tau_1}}_{\tau_1}$ in
(\ref{eq2}) and noting that $1_{[\tau_1<T]}$ is
$F_{\tau_1}$-measurable, it follows that
$$
\begin{array}{lll}Y^1_0&=E[\integ{0}{\tau_1}
\psi_1(s,X_s)ds-\ell_{1 u_{\tau_1}}(\tau_1)\ind_{[\tau _1<T]}]\\ &
+E[\int_{\tau_1}^{\tau_2}\psi_{u_{\tau_1}}(s,X_s)ds.\ind_{[\tau
_1<T]} -\ell_{u_{\tau_1}{u_{\tau_2}}}(\tau_2)\ind_{[\tau_2
<T]}+Y^{u_{\tau_2}}_{\tau_2}\ind_{[\tau_2 <T]}]\\
{}&=E[\integ{0}{\tau_2} \psi_{u_s}(s,X_s)ds-\ell_{1
u_{\tau_1}}(\tau_1)\ind_{[\tau _1<T]}]
-\ell_{u_{\tau_1}{u_{\tau_2}}}(\tau_2)\ind_{[\tau_2
<T]}+Y^{u_{\tau_2}}_{\tau_2}\ind_{[\tau_2 <T]}], \end{array}
$$
since $[\tau_2<T]\subset [\tau_1<T]$.

\noindent Repeating this procedure $n$ times, we obtain \be
\label{eq4} Y^1_0=E[\integ{0}{\tau_n}
\psi_{u_s}(s,X_s)ds-\sum_{j=1}^{n}\ell_{u_{\tau_{j-1}}
u_{\tau_{j}}}(\tau_j)\ind_{[\tau _j<T]}]
+Y^{u_{\tau_n}}_{\tau_n}\ind_{[\tau_n <T]}] \ee But, the strategy
$\delta=(\tau_n)_{n\geq 1}$ is finite, otherwise $Y^1_0$ would be
equal to $-\infty$ contradicting the assumption that the processes
$Y^j$ belong to ${\cal S}^p$. Therefore, taking the limit as
$n\rightarrow \infty$ we obtain $Y^1_0=J(\delta,u)$.

\medskip \noindent To complete the proof it remains to show that
$J(\delta,u)\geq J(\theta,v)$ for any other finite admissible
strategy $(\theta,v):=((\theta_n)_{n\geq 1}, (\zeta_n)_{n\geq 1})$.

\noindent The definition of the Snell envelope yields
$$
\begin{array}{ll}
Y^1_0&\geq E[\int_0^{\theta_1}\psi_1(s,X_s)ds +\max_{j\in {\cal
J}^{-1}}(-\ell_{1j}(\theta_1)+Y^j_{\theta_1})1_{[\theta_1 <T]}]\\
{}&\geq E[\int_0^{\theta_1}\psi_1(s,X_s)ds +(-\ell_{1v_{\theta_1}
}(\theta_1)+Y^{v_{\theta_1}}_{\theta_1})1_{[\theta_1 <T]}].
\end{array}
$$
But, once more using a similar characterization as (\ref{eq3}), we
get
$$
\begin{array}{ll}
Y^{v_{\theta_1}}_{\theta_1}&\geq
E[\int_{\theta_1}^{\theta_2}\psi_{v_{\theta_1}}(s,X_s)ds +\max_{j\in
{\cal
J}^{-{v_{\theta_1}}}}(-\ell_{v_{\theta_1}j}(\theta_2)+Y^j_{\theta_2})1_{[\theta_2
<T]}|\cF_{\theta_1}]\\
{}&\geq E[\int_{\theta_1}^{\theta_2}\psi_{v_{\theta_1}}(s,X_s)ds
+(-\ell_{v_{\theta_1}v_{\theta_2}}(\theta_2)+Y^{v_{\theta_2}}_{\theta_2})1_{[\theta_2
<T]}|\cF_{\theta_1}].
\end{array}
$$
Therefore,
$$
\begin{array}{lll}
Y^1_0&\geq E[\int_0^{\theta_1}\psi_1(s,X_s)ds]
-\ell_{1v_{\theta_1}}(\theta_1)1_{[\theta_1 <T]}] \\ &
+E[1_{[\theta_1
<T]}\int_{\theta_1}^{\theta_2}\psi_{v_{\theta_1}}(s,X_s)ds
-\ell_{v_{\theta_1}v_{\theta_2}}(\theta_2)1_{[\theta_2
<T]}+Y^{v_{\theta_2}}_{\theta_2}1_{[\theta_2 <T]}]\\
{}&=E[\int_0^{\theta_2}\psi_{v_s}(s,X_s)ds
-\ell_{1v_{\theta_1}}(\theta_1)1_{[\theta_1 <T]}
-\ell_{v_{\theta_1}v_{\theta_2}}(\theta_2)1_{[\theta_2
<T]}+Y^{v_{\theta_2}}_{\theta_2}1_{[\theta_2 <T]}].\end{array}
$$
Repeat this argument $n$ times to obtain
$$
Y^1_0\geq E[\int_0^{\theta_n}\psi_{v_s}(s,X_s)ds
-\sum_{j=1}^n\ell_{v_{\theta_{n-1}}v_{\theta_n}}(\theta_n)1_{[\theta_n
<T]} +Y^{v_{\theta_n}}_{\theta_n}1_{[\theta_n <T]}].
$$
Finally, taking the limit as $n\rightarrow \infty$ yields
$$
Y^1_0\geq E[\int_0^{T}\psi_{v_s}(s,X_s)ds -\sum_{j\geq
1}\ell_{v_{\theta_{n-1}}v_{\theta_n}}(\theta_n)1_{[\theta_n <T]}
]=J(\theta,v)$$ since the strategy $(\theta,v)$ is finite. Hence,
the strategy $(\delta,u)$ is optimal. The proof is now complete.
$\Box$

%%%%%%%%%%%%%Existence%%%%%%%%%%%%%%%%%%%%%%%%%%%%

\section{Existence of the processes $(Y^1,\ldots, Y^q)$.}

We will now establish existence of the processes $(Y^1,\ldots,Y^q)$.
They will be obtained as a limit of a sequence of processes
$(Y^{1,n},\dots,Y^{q,n})_{n\geq 0}$ defined recursively as follows.

For $i\in {\cal J}$, set, for every $0\le t\leq T$, \be
\label{eqapp0} Y^{i,0}_t=E[\int_t^T\psi_i(s,X_s)ds|\cF_t], \ee and,
for $n\geq 1$, \be \label{eqappn} Y^{i,n}_t=\mbox{ess sup}_{\tau\geq
t}E[\integ{t}{\tau} \psi_i(s,X_s)ds+\max_{k\in {\cal
J}^{-i}}(-\ell_{ik}(\tau)+Y^{k,n-1}_\tau)\ind_{[\tau<T]}|\cF_t]. \ee

\noindent Set ${\cal D}^{i,n}_t=\{ (\delta,u)=((\tau_n)_{n\geq
1},(\xi_n)_{n\geq 1}) \mbox{ such that } u_0=i,\,\tau_1\geq t\mbox{
and }\tau_{n+1}=T\}$.

\noindent Using the same arguments as the ones of the Verification
Theorem, Theorem \ref{thmverif}, the following characterization of
the processes $Y^{i,n}$ holds true. \be
Y^{i,n}_t=\esssup_{(\delta,u)\in {\cal D}^{i,n}_t}E[\int_t^T
\psi_{u_s} (s,X_s)ds-\sum_{j=1}^n\ell_{u_{\tau_{j-1}}
u_{\tau_{j}}}(\tau_j)\ind_{[\tau _j<T]} |\cF_t]. \ee

\noindent In the next proposition we collect some useful properties
of $Y^{1,n},\ldots,Y^{q,n}$. In particular, we show that, as
$n\to\infty$, the limit processes $\tilde Y^i:=\lim_{n\to
\infty}Y^{i,n}$ exist and are only \cadlag but have the same
Characterization (\ref{eqvt}) as the $Y^i$'s. Thus, the existence
proof of the $Y^i$'s will consist in showing that $\tilde Y^i$'s are
continuous and hence satisfy the Verification Theorem. This will be
done in Theorem 2, below.

\no \bp \label{prop2}
\begin{itemize}
\item[$(i)$] For each $n\geq 0$, the processes $Y^{1,n},\ldots,Y^{q,n}$ are
continuous and belong to ${\cal S}^p$.

\item[$(ii)$] For any $i\in {\cal J}$, the sequence $(Y^{i,n})_{n\geq 0}$
converges increasingly and pointwisely P-$a.s.$ for any $0\le t\leq
T$ and in ${\cal M}^{p,1}$ to \cadlag processes $\tilde{Y}^i$.
Moreover, these limit processes satisfy
\begin{itemize}
\item[$(a)$]
\begin{equation}
E[\sup_{0\le t\leq T}|\tilde{Y}^i_t|^p]<\infty,\,\,\quad  i\in {\cal
J}.
\end{equation}
\item[$(b)$] For any $0\le t\leq T$,
\begin{equation}
\label{eq26} \tilde{Y}^{i}_t=\mbox{ess sup}_{\tau\geq
t}E[\integ{t}{\tau} \psi_i(s,X_s)ds+\max_{k\in {\cal
J}^{-i}}(-\ell_{ik}(\tau)+\tilde{Y}^{k}_\tau)\ind_{[\tau<T]}|F_t].
\end{equation}
\end{itemize}
\end{itemize}
\ep

\no {\it Proof}. $(i)$ Let us show by induction that, for any $n\geq
0$, $Y^{i,n}_T=0$ and $Y^{i,n}\in {\cal S}^p$, for every $i\in {\cal
J}$.

\noindent For $n=0$ the property holds true since the process
$(\psi_i(s,X_s))_{0\le s\leq T}$ belongs to ${\cal S}^p$. Suppose
now that the property is satisfied for some $n$. By Proposition
\ref{thmsnell}, for every $i\in{\cal J}$ and up to a term,
$Y^{i,n+1}$ is the Snell envelope of the process $\left(\integ{0}{t}
\psi_i(s,X_s)ds+\max_{k\in {\cal
J}^{-i}}(-\ell_{ik}(t)+Y^{k,n}_t)\ind_{[t<T]}\right)_{0\le t\leq T}$
and verifies $Y^{i,n+1}_T=0$. Since $\max_{k\in {\cal
J}^{-i}}(-\ell_{ik}(t)+Y^{k,n}_t)_{\big|_{t=T}}<0$, this process is
continuous on $[0,T)$ and have a positive jump at $T$, $Y^{i,n+1}$
is continuous and belongs to ${\cal S}^p$. This shows that, for
every $i\in {\cal J}$, $Y^{i,n}_T=0$ and $Y^{i,n}\in {\cal S}^p$ for
any $n\geq 0$.

\medskip

\noindent $(ii)$ We show by induction on $n\geq 0$, that for each
$i\in {\cal J}$,
$$
Y^{i,n}\leq Y^{i,n+1}\leq
E[\int_t^T\max_{i=1,\ldots,q}|\psi_i(s,X_s)|ds|\cF_t].
$$
For $n=0$ the property is obviously true, since it is enough to take
$\tau=T$ in the definition of $Y^{i,1}$ to obtain that $Y^{i,1}\geq
Y^{i,0}$. On the other hand taking into account that $\ell_{ij}\ge
\gamma> 0$ we have \be \label{y1}
\begin{array}{ll}
Y^{i,1}_t&=\mbox{ess sup}_{\tau\geq t}E[\integ{t}{\tau}
\psi_i(s,X_s)ds+\max_{k\in {\cal
J}^{-i}}(-\ell_{ik}(\tau)+Y^{k,0}_\tau)\ind_{[\tau<T]}|\cF_t]\\{}&\leq
\mbox{ess sup}_{\tau\geq t}E[\integ{t}{\tau} \psi_i(s,X_s)ds+
E[\int_\tau^T\max_{i=1,\ldots,q}|\psi_i(s,X_s)|ds|F_\tau]|\cF_t]
\\ &\leq
E[\int_t^T\max_{i=1,\ldots,q}|\psi_i(s,X_s)|ds|\cF_t].
\end{array}
\ee

\noindent Suppose now that, for some $n$, we have
$$Y^{i,n}\leq Y^{i,n+1}\leq
E[\int_t^T\max_{i=1,\ldots,q}|\psi_i(s,X_s)|ds|\cF_t],\,\,\quad i\in
{\cal J}.
$$
Replace $Y^{i,n+1}$ by $Y^{i,n}$ in the definition of $Y^{i,n+2}$,
to obtain that $Y^{i,n+2}\geq Y^{i,n+1}$.

\noindent Finally, as is the case for $Y^{i,1}$ in (\ref{y1}), we
also have
$$
Y^{i,n+2}_t\leq
E[\int_t^T\max_{i=1,\ldots,q}|\psi_i(s,X_s)|ds|\cF_t],\,\,\quad 0\le
t\leq T.
$$
Therefore, for every $i\in {\cal J}$, the sequence $(Y^{i,n})_{n\geq
0}$ is increasing in $n$ and satisfies \be \label{inq} Y^{i,n}_t\leq
E[\int_t^T\max_{i=1,\ldots,q}|\psi_i(s,X_s)|ds|\cF_t],\,\,\quad
0\le, t\leq T. \ee Therefore, it converges to some limit
$\tilde{Y}^i_t:=\lim_{n\rightarrow \infty}Y^{i,n}_t$ that satisfies
$$
Y^{i,0}_t\leq \tilde{Y}^i_t\leq
E[\int_t^T\max_{i=1,\ldots,q}|\psi_i(s,X_s)|ds|\cF_t], \,\quad 0\le
t\leq T.
$$
Now, using the smoothness properties of $\psi_i$,  Doob's Maximal
Inequality yields that, for each $i\in {\cal J}$,
$$
E[\sup_{0\le t\leq T}|\tilde{Y}^i_t|^p]<\infty.
$$
By the Lebesgue Dominated Convergence Theorem, the sequence
$(Y^{i,n})_{n\geq 0}$ also converges to $\tilde{Y}^i$ in ${\cal
M}^{p,1}$.

\noindent Let us now show that $\tilde{Y}^i$ is \cadlag. We note
that, for each $n\geq 1$ and $i\in {\cal J}$, the process
$\left(Y^{i,n}_t+\int_0^t\psi_i(s,X_s)ds\right)_{0\le t\leq T}$ is a
continuous supermartingale, since, by Eq.\@ (\ref{eqappn}), it is
the Snell envelope of the continuous process $\left(\integ{0}{t}
\psi_i(s,X_s)ds+\max_{k\in {\cal
J}^{-i}}(-\ell_{ik}(\tau)+Y^{k,n-1}_t)\ind_{[t<T]}\right)_{0\le
t\leq T}$. Hence, its limit process
$\left(\tilde{Y}^{i}_t+\int_0^t\psi_i(s,X_s)ds\right)_{0\le t \leq
T}$ is \cadlag, as a limit of increasing sequence of continuous
supermartingales. Therefore, $\tilde{Y}^i$ is \cadlag.

\noindent Finally, the \cadlag processes
$\tilde{Y}^1,\ldots,\tilde{Y}^q$ satisfy Eq.\@ (\ref{eq26}), since
they are limits of the increasing sequence of processes $Y^{i,n}$,
$i\in {\cal J}$, that satisfy (\ref{eqappn}). We use Proposition
\ref{thmsnell}, $(v)$ to conclude. $\Box$

\medskip

We will now prove that the processes $\tilde{Y}^1,\ldots,\tilde{Y}^q$
are continuous and satisfy the Verification Theorem, Theorem
\ref{thmverif}.

\begin{th} \label{thm3}The limit processes
$\tilde{Y}^1,\ldots,\tilde{Y}^q$ satisfy the Verification Theorem.
\end{th}

\no {\it Proof}. Recall from Proposition \ref{prop2} that the
processes $\tilde{Y}^1,\ldots,\tilde{Y}^q$ are \cadlag, uniformly
$L^p$-integrable and satisfy (\ref{eq26}). It remains to prove that
they are continuous.

\noindent Indeed, note that, for $i\in {\cal J}$, the process
$\left(\tilde{Y}^i_t+\int_0^t\psi_i(s,X_s)ds\right)_{0\le t\leq T}$
is the Snell envelope of
$$
\left(\int_0^t\psi_i(s,X_s)ds+ \max_{k\in {\cal
J}^{-i}}(-\ell_{ik}(t)+\tilde{Y}^{k}_t)\ind_{[t<T]}\right)_{0\le
t\leq T}.
$$
Therefore, thanks to the Doob-Meyer decomposition of the Snell
Envelope of processes (Proposition \ref{thmsnell}-$(iii)$), there
exist continuous martingales $(M^i_t)_{t\leq T}$ and continuous,
resp. purely discontinuous, nondecreasing processes $(A^i_t)_{t\leq
T}$, resp. $(B^i_t)_{t\leq T}$, such that, for each $i\in {\cal J}$,
and $0\le t\leq T$,
$$
\begin{array}{l}
\int_0^t\psi_i(s,X_s)ds+\tilde{Y}^i_t=M^i_t-A^i_t-B^i_t \,\,\quad
(A^i_0=B^i_0=0).
\end{array}
$$
Moreover, the following properties for the jumps of $B^i$, $i\in
{\cal J}$ hold. When there is a jump of $B^i$ at $t$, there is a
jump, at the same time $t$, of the process $(\max_{k\in {\cal
J}^{-i}}(-\ell_{ik}(t)+\tilde{Y}^{k}_t))_{t\leq T}$. Since
$\ell_{ij}$ are continuous, there is $j\in {\cal J}^{-i}$ such that
$\Delta_t\tilde{Y}^j=-\Delta_tB^j<0$ and
$\tilde{Y}^i_{t-}=-\ell_{ij}(t)+\tilde{Y}^j_{t-}$. Suppose now there
is an index $i_1\in {\cal J}$ for which there exists $t\in [0,T]$
such that $\Delta_tB^{i_1}>0$. This implies that there exists
another index $i_2\in {\cal J}^{-i_1}$ such that $\Delta_t
B^{i_2}>0$ and
$\tilde{Y}^{i_1}_{t-}=-\ell_{i_1i_2}(t)+\tilde{Y}^{i_2}_{t-}$. But,
given $i_2$, there exists an index $i_3\in {\cal J}^{-i_2}$ such
that $\Delta_t B^{i_3}>0$ and
$\tilde{Y}^{i_2}_{t-}=-\ell_{i_2i_3}(t)+\tilde{Y}^{i_3}_{t-}$.
Repeating this argument many times, we get a sequence of indices
$i_1,\ldots,i_j,\ldots \in {\cal J}$ that have the property that
$i_{k}\in {\cal J}^{-i_{k-1}}$, $\Delta_t B^{i_k}>0$ and
$\tilde{Y}^{i_{k-1}}_{t-}=-\ell_{i_{k-1}i_k}(t)+\tilde{Y}^{i_k}_{t-}$.

\noindent Since ${\cal J}$ is finite then there exist two indices
$m<r$ such that $i_m=i_r$ and $i_m, i_{m+1}, ..., i_{r-1}$ are
mutually different. It follows that:
$$
\tilde{Y}^{i_m}_{t-}=-\ell_{i_mi_{m+1}}(t)+\tilde{Y}^{i_{m+1}}_{t-}
=-\ell_{i_mi_{m+1}}(t)
-\ell_{i_{m+1}i_{m+2}}(t)+\tilde{Y}^{i_{m+2}}_{t-}=\cdots=
-\ell_{i_mi_{m+1}}(t)-\cdots-\ell_{i_{r-1}i_r}(t)+\tilde{Y}^{i_{r}}_{t-}.
$$
As $i_m=i_r$ we get
$$
-\ell_{i_mi_{m+1}}(t)-\cdots-\ell_{i_{r-1}i_r}(t)=0
$$
which is impossible since for any $i\neq j$, all $0\le t\leq T$,
$\ell_{ij}(t)\geq \gamma >0$. Therefore, there is no $i\in {\cal J}$
for which there is a $t\in [0,T]$ such that $\Delta_tB^i>0$. This
means that $B^i\equiv 0$ and the processes
$\tilde{Y}^1,\ldots,\tilde{Y}^q$ are continuous. Since they satisfy
(\ref{eq26}), then, by uniqueness, $Y^i=\tilde{Y}^i$, for any $i\in
{\cal J}$. Thus, the Verification Theorem \ref{thmverif} is
satisfied by $Y^1,\ldots,Y^q$. $\Box$
\medskip

We end this section by the following convergence result of the
sequences $(Y^{i,n})_{n\geq 0}$ to $Y^i$'s.
\begin{pro}\label{propconv} It holds true that, for any $i\in {\cal J}$,
$$
E[\sup_{s\leq T}|Y^{i,n}_s-Y^i_s|^p]\to 0 \quad\mbox{as}\quad n\to
+\infty.
$$
\end{pro}
{\it Proof}. By Proposition \ref{prop2}, we know that $P$-a.s., for
any $n\geq 1$, the function $t\mapsto Y^{i,n}_t(\omega)$ is
continuous and for any $0\le t\leq T$ the sequence
$(Y^{i,n}_t(\omega))_{n\geq 1}$ converges increasingly to
$Y^i_t(\omega)$. As the function $t\mapsto Y^{i}_t(\omega)$ is
continuous then thanks to Dini's Theorem it holds true that:
$$
P-\mbox{a.s.}\,\, \lim_{n\to\infty}\sup_{0\le t\leq
T}|Y^{i,n}_t(\omega)-Y^i_t(\omega)|=0.
$$ The result now follows from
the Lebesgue Dominated Convergence Theorem. $\Box$

%%%%%%%%%%%%Connection with QVI%%%%%%%%%%%%%%%

\section{Connection with systems of variational inequalities} When the
underlying market price process $X$ is Markov diffusion and the
switching costs are of the form $\ell_{ij}(t,X_t)$, the classical
methods of solving impulse problems (cf. Brekke and \O ksendal
(1994), Guo and Pham (2005)) formulate a Verification Theorem
suggesting that the value function of our optimal switching problem
is the unique viscosity solution the following system of
quasi-variational inequalities (QVI) with inter-connected obstacles
\be \label{qv1}\left\{
\begin{array}{ll}
\min\{\phi_i(t,x)-\max_{j\in {\cal
J}^{-i}}(-\ell_{ij}(t,x)+\phi_j(t,x)),-\partial_t\phi_i(t,x)-A\phi_i(t,x)-\psi_i(t,x)\}=0,\\
\phi_i(T,x)=0, \quad\quad i\in {\cal J},
\end{array}\right.\ee
where $A$ is the infinitesimal generator of the driving process $X$.

\medskip\noindent
However, besides the technical difficulties to establish existence
of a smooth solution, existence and uniqueness of a viscosity
solution for such systems still remains open for most of the models
discussed in the literature (See Carmona and Ludkovski (2006) for a
detailed discussion).

\medskip
By means of yet another characterization of the Snell envelope in
terms of systems of reflected Backward SDEs, due to El Karoui et al.
(1997-1)(Theorems 7.1 and 8.5), we are able to show that the vector
of value processes $(Y^1,\ldots,Y^q)$ of our optimal problem is a
viscosity solution of the system (\ref{qv1}), when the switching
cost functions $\ell_{ij}$ are only deterministic functions of {\it
the time variable}. An example of such a family of switching costs
is
$$
\ell_{ij}(t)=e^{-rt}a_{ij}, \quad
$$
where, $a_{ij}$ are constant costs and $r>0$ is some discounting
rate.

We show that under mild assumptions on the coefficients
$\psi_i(t,x)$ and $\ell_{ij}(t)$,
$$
Y_t^i=v^i(t,X_t),\,\,\quad 0\le t\leq T,\,\, \, i\in {\cal J},
$$
where the deterministic functions $v^1(t,x),\ldots, v^q(t,x)$ are
viscosity solutions of the following system of QVI with
inter-connected obstacles \be \label{sysvi}\left\{
\begin{array}{l}
\min\{v_i(t,x)-\max_{j\in {\cal
J}^{-i}}(-\ell_{ij}(t)+v_j(t,x)),-\partial_tv_i(t,x)-Av_i(t,x)-\psi_i(t,x)\}=0,\\
v_i(T,x)=0, \quad\quad i\in {\cal J}.
\end{array}\right.
\ee

For $(t,x)\in [0,T]\times \R^k$, let $(X^{tx}_s)_{s\leq T}$ be the
solution of the following It\^o diffusion:
\begin{equation}\label{sde}
dX^{tx}_s=b(s,X_s^{tx})ds+\sigma(s,X_s^{tx})dB_s,\quad \,t\leq s\leq
T;\,\,\quad\quad X_s^{tx}=x\,\,\, \mbox{ for }s\leq t,
\end{equation}
where, the functions $b$ and $\sigma$, with appropriate dimensions,
satisfy the following standard conditions:

There exists a constant $C\geq 0$ such that \be
\label{regbs}|b(t,x)|+ |\sigma(t,x)|\leq C(1+|x|) \quad \mbox{ and }
\quad |\sigma(t,x)-\sigma(t,x')|+|b(t,x)-b(t,x')|\leq C|x-x'|\ee for
any $t\in [0, T]$ and $x, x'\in \R^k$.

\noindent These properties of $\sigma$ and $b$ imply in particular
that the process $X^{tx}:=(X^{tx}_s)_{0\le s\leq T}$, solution of
(\ref{sde}), exists and is unique. Its infinitesimal generator $A$
is given by
\begin{equation}\label{generator}
A=\frac{1}{2}\sum_{i,j=1}^d(\sigma
.\sigma^*)_{ij}(t,x)D_{ij}+\sum_{i=1}^d b_i(t,x)D_i.
\end{equation}
Moreover, the following estimates hold true (see e.g. Revuz and Yor
(1991) for more details).
\medskip
\bp \label{estimx} The process $X^{tx}$ satisfies the following
estimates:
\begin{itemize}
\item [$(i)$] For any $\theta\geq 2$, there exists a constant $C$ such that
\be \label{estimx0}E[\sup_{0\le s\leq T}|X^{tx}_s|^\theta]\leq
C(1+|x|^\theta). \ee
\item[$(ii)$] There exists a constant $C$ such that for any $t,t'\in
[0,T]$ and $x,x'\in \R^k$, \be \label{estimx1} E[\sup_{0\le s\leq
T}|X^{tx}_s-X^{t'x'}_s|^2]\leq
C(1+|x|^2)(|x-x'|^2+|t-t'|).\,\,\,\Box \ee
\end{itemize}
\ep

Let us now introduce the following assumption on the payoff rates
$\psi_i$ and the switching cost functions $\ell_{ij}$:
\medskip

\noindent {\bf Assumption [H]}.
\begin{itemize}
\item [(H1)] The running costs
$\psi_i$, $i=1,\ldots,q$, (of Subsection 2.1) are jointly continuous
and are of polynomial growth, i.e., there exist some positive
constants $C$ and $\delta$ such that for each $i\in {\cal J}$,
$$
|\psi_i(t,x)|\leq C(1+|x|^{\delta}),\,\,\quad (t,x)\in [0,T]\times
\R^k.
$$

\item [(H1)] For any $i,j\in {\cal J}$, the switching costs
$\ell_{ij}$ are deterministic functions of $t$ and continuous and
there exists a real constant $\gamma>0$ such for any $0\le t\leq T$,
$\min\{\ell_{ij}(t), i,j\in {\cal J}, \,i\neq j\}\geq \gamma$.
\end{itemize}

\medskip\noindent Taking into account Proposition \ref{estimx}, the processes
$(\psi_i(s,X^{tx}_s)_{0\le s\leq T})_{i=1,q}$ belong to ${\cal M}^{2,1}$. A
condition we will need to establish a characterization of the value
processes of our optimal problem with a class of reflected backward
SDEs. Note that the required polynomial growth condition on the
$\psi_i$'s is not contradictory with the condition listed in
Assumptions 2.1 (ii), since the process $X^{tx}$ has finite moments
of all orders (see also Remark 1).
\medskip

Recall the notion of viscosity solution of the system (\ref{sysvi}).
\begin{axiom}
Let $(v_1,\ldots,v_q)$ be a vector of continuous functions on
$[0,T]\times \R^k$ with values in $\R^q$ and such that
$(v_1,\ldots,v_q)(T,x)=0$ for any $x\in \R^k$. The vector
$(v_1,\ldots,v_q)$ is called:
\begin{itemize}
\item [$(i)$] A viscosity supersolution of the system $(\ref{sysvi})$ if for
any $(t_0,x_0)\in [0,T]\times \R^k$ and any $q$-tuplet functions
$(\varphi_1,\dots,\varphi_q)\in (C^{1,2}([0,T]\times \R^k))^q$ such
that $(\varphi_1,\dots,\varphi_q)(t_0,x_0)=(v_1,\dots,v_q)(t_0,x_0)$
and for any $i\in {\cal J}$, $(t_0,x_0)$ is a maximum of $\varphi_i
-v_i$ then we have: for any $i\in {\cal J}$, \be
\min\{v_i(t_0,x_0)-\max_{j\in {\cal
J}^{-i}}(-\ell_{ij}(t_0)+v_j(t_0,x_0)),-\partial_t\varphi_i(t_0,x_0)-A\varphi_i(t_0,x_0)-\psi_i(t_0,x_0)\}\geq
0. \ee
\item [$(ii)$] A viscosity subsolution of the system $(\ref{sysvi})$ if for
any $(t_0,x_0)\in [0,T]\times \R^k$ and any $q$-tuplet functions
$(\varphi_1,\dots,\varphi_q)\in (C^{1,2}([0,T]\times \R^k))^q$ such
that $(\varphi_1,\dots,\varphi_q)(t_0,x_0)=(v_1,\dots,v_q)(t_0,x_0)$
and for any $i\in {\cal J}$, $(t_0,x_0)$ is a minimum of $\varphi_i
-v_i$ then we have: for any $i\in {\cal J}$, \be
\min\{v_i(t_0,x_0)-\max_{j\in {\cal
J}^{-i}}(-\ell_{ij}(t_0)+v_j(t_0,x_0)),-\partial_t\varphi_i(t_0,x_0)-A\varphi_i(t_0,x_0)-\psi_i(t_0,x_0)\}\leq
0. \ee
\item [$(iii)$] The vector of function $(v_1,\ldots,v_q)$ is a viscosity
solution of the system $(\ref{sysvi})$ if it is both a viscosity
supersolution and subsolution.
\end{itemize}
\end{axiom}

Let now $(Y^{1,tx}_s,\ldots,Y^{q,tx}_s)_{0\le s\leq T}$ be the
vector of value processes which satisfies the Verification Theorem
\ref{thmverif} associated with $(\psi_i(s,X^{tx}_s))_{s\leq T}$ and
$\ell^{ij}(t)$. The vector $(Y^{1,tx},\ldots, Y^{q,tx})$ exists
through Theorem \ref{thm3} combined with the estimates of $X^{tx}$
of Proposition \ref{estimx} and Assumptions [H].
\medskip

The following theorem is the main result of this section.

\beth Under Assumption $[\bf H]$, there exist $q$ deterministic
functions $v^1(t,x),\dots,v^q(t,x)$ defined on $[0,T]\times \R^k$
and $\R$-valued such that:

\begin{itemize}
\item [$(i)$] $v^1,\ldots,v^q$ are continuous in $(t,x)$,
are of polynomial growth and satisfy, for each $t\in [0, T]$, and
for every $s\in [t,T]$,
$$
Y^{i,tx}_s=v^i(s,X_s^{tx}),\,\,\mbox{ for every } i\in {\cal J}.$$
\item [$(ii)$] The vector of functions $(v^1,\ldots,v^q)$ is a viscosity
solution for the system of variational inequalities (\ref{sysvi}).
\end{itemize}
\eeth

\noindent {\it Proof}. The proof is obtained through the three
following steps.

\noindent{\bf Step 1}. {\it An approximation scheme}
\medskip

\noindent For $n\geq 0$, let $(Y^{1,n,tx}_s)_{ 0\le s\leq T},\ldots,
(Y^{q,n,tx}_s)_{  0\le s\leq T}$ be the continuous processes defined
recursively by Eqs.\@ (\ref{eqapp0})-(\ref{eqappn}). Using
Assumption [H1], the estimates (\ref{estimx0}) for $X^{tx}$ and
Proposition \ref{prop2}, the processes $Y^{1,n,tx}\dots, Y^{q,n,tx}$
belong to ${\cal S}^2$. Therefore, using a result by El Karoui et
al. ((1997-1), Theorem 7.1) which characterizes a Snell envelope as
a solution for a one barrier reflected BSDE, for any $n\geq 1$ and
$i\in {\cal J}$, there exists a pair of ${\cal F}_t$-adapted
processes $(Z^{i,n,tx},K^{i,n,tx})$ with value in $R^d\times R^+$
such that:
\begin{equation}\label{yn2} \left\{
\begin{array}{l}
Y^{i,n,tx},\, K^{i,n,tx}\in {\cal S}^2 \, \mbox{ and }\,
Z^{i,n,tx}\in {\cal M}^{2,d};\,K^{i,n,tx}
\mbox{ is  nondecreasing and }K^{i,n,tx}_0=0,\\
Y^{i,n,tx}_s=\integ{s}{T}\psi_i(u,X_u^{tx})du-\integ{s}{T}Z^{i,n,tx}_udB_u+K_T^{i,n,tx}-K^{i,n,tx}_s,
\,\, \mbox{for all} \,\,\, 0 \le s\le T,\\ Y^{i,n,tx}_s\geq
\max_{j\in {\cal J}^{-i}}\{-\ell_{ij}(s)+Y^{j,n-1,tx}_s\},\,\,\,\,
\mbox{for all} \,\,\, 0 \le s\le
T,\\
\integ{0}{T}(Y^{i,n,tx}_u- \max_{j\in {\cal
J}^{-i}}\{-\ell_{ij}(u)+Y^{j,n-1,tx}_u\})dK^{i,n}_u=0.
\end{array}
\right. \end{equation}

Thanks to Theorem 8.5 in  El Karoui et al.\@ (1997-1) related to the
representation of solutions of reflected backward SDEs, there exist
deterministic functions $v^{1,0}, \dots, v^{q,0}$ defined on
$[0,T]\times R^k$, continuous and with polynomial growth such that
for every $(t,x)\in [0,T]\times R^k$ and every $i\in {\cal J}$,
$$
Y^{i,0,tx}_s=v^{i,0}(s,X^{tx}_s),\,\,\quad t\le s\le T.
$$

Using an induction argument, and applying Theorem 8.5 in  El Karoui
et al.\@ (1997-1) at each step, yields the existence of
deterministic functions $v^{1,n}, \ldots, v^{q,n}$ defined on
$[0,T]\times R^k$, that are continuous and with polynomial growth
such that, for every $(t,x)\in [0,T]\times R^k$ and every $i\in
{\cal J}$,
$$
Y^{i,n,tx}_s=v^{i,n}(s,X^{tx}_s),\,\,\quad t\le s\le T.
$$
Since the sequences of processes $(Y^{i,n,tx})_{n\geq 0}$ is
nondecreasing in $n$, then for any $i\in {\cal J}$, the sequences of
deterministic functions $(v^{i,n})_{n\geq 0}$ is also nondecreasing.

Moreover, we have \be\label{v}
\begin{array}{lll}
v^{i,n}(t,x)&\leq Y^{t,x}_t\leq
E[\int_t^T\max_{i=1,\ldots,q}|\psi_i(s,X^{tx}_s)|ds|\cF_t]\\
{}&\leq E[\int_t^T\max_{i=1,\ldots,q}|\psi_i(s,X^{tx}_s)|ds],
\end{array}
\ee where, the last inequality is obtained after taking
expectations, since $v^{i,n}(t,x)$ is a deterministic function. It
follows that for any $i\in {\cal J}$, the sequence $(v^{i,n})_{n\geq
0}$ converges pointwisely to a deterministic function $v^i$ and the
last inequality in Eq. (\ref{v}) implies that $v^i$ is of polynomial
growth through $\psi_i$ and the estimates (\ref{estimx0}) for
$X^{tx}$. Furthermore, for any $(t,x)\in [0,T]\times R^k$ we have
\be \label{Y=v} Y^{i,tx}_s=v^i(s,X^{tx}_s),\,\,\quad t\le s\le T.
\ee
\medskip

\noindent {\bf Step 2}. {\it $L^2(P)$-continuity of the value
functions $(t,x)\longrightarrow Y^{i,tx}$}.

\medskip\noindent
Let $(t,x)$ and $(t',x')$ be elements of  $[0,T]\times \R^k$. Using
the representation (\ref{eqyn}) we will show that
$$
E[\sup_{0\le s\leq T}|Y^{i,t'x'}_s-Y^{i,tx}_s|^2]\rightarrow 0\mbox{
as }(t',x')\rightarrow (t,x) \mbox{ for any }i\in {\cal J}.
$$
Indeed, recall that, by (\ref{eqyn}), we have, for any $i\in {\cal
J}$ and $s\in [0,T]$
$$
Y^{i,tx}_s=\esssup_{(\delta,u)\in {\cal
D}_s}E[\int_s^T\psi_{u_s}(s,X^{tx}_s)ds- \sum_{j\geq
1}\ell_{u_{\tau_{j-1}} u_{\tau_{j}}}(\tau_j)\ind_{[\tau _j<T]}
|\cF_s],
$$
where, ${\cal D}_s$ is the set of finite strategies such that
$\tau_1\geq s$, $P-a.s.$

\medskip\noindent Therefore,
$$
\begin{array}{ll}
|Y^{i,tx}_s-Y^{i,t'x'}_s|&\leq \esssup_{(\delta,u) \in {\cal
D}_s}E[\integ{s}{T}|\psi_{u_r}(r,X_r^{tx})-\psi_{u_r}(r,X_r^{t'x'},u_r)|dr|{\cal
F}_s]\\
{}&\leq
E[\integ{0}{T}\{\sum_{i=1}^q|\psi_i(r,X_r^{tx})-\psi_i(r,X_r^{t'x'})|\}ds|{\cal
F}_s].
\end{array}
$$
Now, using Doob's Maximal Inequality (see e.g. \cite{[RY]}) and
taking expectation, there exists of a constant $C\geq 0$ such that:
\be \label{eqinter}E[\sup_{0\le s\leq
T}|Y^{i,tx}_s-Y^{i,t'x'}_s|^2]\leq
CE[\integ{0}{T}\{\sum_{i=1}^q|\psi_i(r,X_r^{tx})-\psi_i(r,X_r^{t'x'})|\}^2ds].
\ee But, the right-hand side of this last inequality converges to
$0$ as $(t',x')$ tends to $(t,x)$. Indeed, for any $ \varpi >0$ it
holds true that:
$$
\begin{array}{lll}
E[\integ{0}{T}\{\sum_{i=1}^q|\psi_i(r,X_r^{tx})-\psi_i(r,X_r^{t'x'})|\}^2ds]
& \leq
E[\integ{0}{T}\{\sum_{i=1}^q|\psi_i(r,X_r^{tx})-\psi_i(r,X_r^{t'x'})|\}2\ind_{[
|X_r^{tx}|+|X_r^{t'x'}|\leq \varpi]}ds]\\ & +
E[\integ{0}{T}\{\sum_{i=1}^q|\psi_i(r,X_r^{tx})-\psi_i(r,X_r^{t'x'})|\}2\ind_{[
|X_r^{tx}|+|X_r^{t'x'}|> \varpi]}ds].\end{array}
$$
By the Lebesgue Dominated Convergence Theorem, the continuity of
$\psi_i$ and  Estimates (\ref{estimx1}), the first term of the
right-hand side of this inequality converges to $0$ as $(t',x')$
tends to $(t,x)$.

\medskip\noindent
The second term satisfies:
$$\begin{array}{l}
E[\integ{0}{T}\{\sum_{i=1,\ldots,q}|\psi_i(r,X_r^{tx})-\psi_i(r,X_r^{t'x'})|\}2\ind_{[
|X_r^{tx}|+|X_r^{t'x'}|> \varpi]}ds]\\ \qquad \qquad\qquad \qquad
\leq
\{E[\integ{0}{T}\{\sum_{i=1,\ldots,q}|\psi_i(r,X_r^{tx})-\psi_i(r,X_r^{t'x'})|\}4]\}^{\frac{1}{2}}
\{E[\integ{0}{T}\ind_{[ |X_r^{tx}|+|X_r^{t'x'}|>
\varpi]}ds]\}^{\frac{1}{2}}\\
\qquad \qquad\qquad \qquad\leq
\{E[\integ{0}{T}\{\sum_{i=1,\ldots,q}|\psi_i(r,X_r^{tx})-\psi_i(r,X_r^{t'x'})|\}4]\}^{\frac{1}{2}}
\{\varpi^{-1}E[\integ{0}{T}
(|X_r^{tx}|+|X_r^{t'x'}|)ds]\}^{\frac{1}{2}}.
\end{array}
$$
Using Estimates (\ref{estimx0}) and the polynomial growth of
$\psi_i$, it follows that, when $(t',x')$ tends to $(t,x)$, the
supremum limit of the right-hand side of the last inequality is
smaller than $ \varpi^{-\frac{1}{2}}C_{tx}$ where $C_{tx}$ is a
constant. As $\varpi$ is whatever then going back to (\ref{eqinter})
and taking the limit to obtain, for any $i\in {\cal J}$,
$$
E[\sup_{0\le s\leq T}|Y^{i,tx}_s-Y^{i,t'x'}_s|^2]\rightarrow 0
\,\,\quad \mbox{ as }\,\,\, (t',x')\rightarrow (t,x).$$

\noindent {\bf Step 3}. {\it the functions $v^1,\ldots, v^q$ are
continuous in $(t,x)$ and the vector of functions $(v^1,\ldots,v^q)$
is a viscosity solution of the system of variational inequalities
(\ref{sysvi}).}

\medskip
\noindent Thanks to the result obtained in Step 2, for any $i\in
{\cal J}$, the function $(s,t,x)\mapsto Y^{i,tx}_s$ is continuous
from $[0,T]^2\times \R^k$ into $L^2(P)$. Indeed, this follows from the
fact that
$$
|Y^{i,t'x'}_{s'}-Y^{i,tx}_s|\leq
|Y^{i,t'x'}_{s'}-Y^{i,tx}_{s'}|+|Y^{i,tx}_{s'}-Y^{i,tx}_{s}|\leq
\sup_{s\leq
T}(|Y^{i,t'x'}_{s}-Y^{i,tx}_{s}|)+|Y^{i,tx}_{s'}-Y^{i,tx}_{s}|.
$$
Therefore, the function $(t,t,x)\mapsto Y^{i,tx}_t$ is also
continuous. But, the result obtained in Step 1, implies that
$Y^{i,tx}_t$ is deterministic and is equal to $v^i(t,x)$. Hence, the
function $v^i$ is continuous in $(t,x)$. The deterministic functions
$v^i$, $i\in {\cal J}$, being continuous and of polynomial growth,
by Theorem 8.5 in El-Karoui et al.\@ (1997-1), these functions are
viscosity solutions for the system (\ref{sysvi}). $\Box$
%%%%%%%%%%%%Connection with QVI%%%%%%%%%%%%%%%
%%%%%%%%%%%%%%%%%%%%%%%%%%%%Numerics of RBSDEs-New added by Said July 1,
2007%%%%%%%%%%%%%%%%%%%%%

\section{Simulating the value-processes $(Y^1,...,Y^q)$}
An important issues in the optimal multiple switching
problem is to provide efficient algorithms to simulate of the value-processes
$(Y^1,...,Y^q)$ solution
of the Verification Theorem \ref{thmverif}.
In this section we comment on this by providing yet another approximation scheme
of the
value-processes $(Y^1,...,Y^q)$ by exploiting their representation as solution
for a system of BSDE with one reflecting barrier. Thanks to a result in
El-Karoui et al. ((1997-1), Theorem 7.1) which
characterizes a Snell envelope of a process which belongs to ${\cal
S}^2$ as a solution for a BSDE with one reflecting barrier, the
vector $(Y^1,...,Y^q)$ is the solution of the following system
of reflected BSDEs:

\noindent
For any $i\in {\cal J}$, there exists a pair of
${\cal F}_t$-adapted processes $(Z^{i},K^{i})$ with value in
$\R^d\times \R^+$ such that:
\begin{equation}\label{rBSDE1} \left\{
\begin{array}{l}
Y^{i},\, K^{i}\in {\cal S}^2 \, \mbox{ and }\, Z^{i}\in {\cal
M}^{2,d}\,\,;\,K^{i} \mbox{ is continuous
nondecreasing and }K^{i}_0=0,\\
Y^{i}_s=\integ{s}{T}\psi_i(u,X_u)du-\integ{s}{T}Z^{i}_udB_u+K_T^{i}-K^{i}_s,
\,\, \mbox{for all} \,\,\, 0 \le s\le T,\\ Y^{i}_s\geq \max_{j\in
{\cal J}^{-i}}\{-\ell_{ij}(s)+Y^{j}_s\},\,\,\,\, \mbox{for all}
\,\,\, 0 \le s\le
T,\\
\integ{0}{T}(Y^{i}_u- \max_{j\in {\cal
J}^{-i}}\{-\ell_{ij}(u)+Y^{j}_u\})dK^{i}_u=0.
\end{array}
\right. \end{equation}
Note that, when $X\equiv X^{tx}$,
taking the limit in (\ref{yn2}), we obtain the solution of
the system (\ref{rBSDE1}).

It is now well known that the solution of a reflected BSDE can be
approximated, in using a penalization scheme, by solutions of
standard BSDEs (see El-Karoui et al. (1997-1) for more
details). Indeed, for $n\geq 0$, consider the following sequence of SDEs
\begin{equation} \label{rBSDE2}
Y^{i,n}_t =
\int_t^T \psi_i(s,X_s)ds + n \int_t^T (L^{i,n}_s - Y^{i,n}_s)^+ ds -
\int_t^T Z^{i,n}_s dB_s, \,\,\,  i \in {\cal J}, \,\,\, t \in [0,T],
\end{equation}
where, for every $i \in {\cal J}$,
$$
L^{i,n}_t = \max_{k\in {\cal
J}^{-i}}(-\ell_{ik}(t)+Y^{k,n}_t),\,\,\, t \in [0,T].
$$
Now, if we define the generator $f_n=(f^1_n,\ldots,f^q_n):\,\, [0,T] \times \R^q
\to \R^q$ by
$$
f^i_n(s,(y_1,...,y_q)) = \psi_i(s,X_s) +
n(\max_{k\in {\cal J}^{-i}}(-\ell_{ik}(s)+y_k) - y_i)^+,\,\,\,\,  i \in {\cal
J},
$$
the $\R^q$-valued process $Y^n=(Y^{1,n},\ldots,Y^{q,n})$ satisfies the
following BSDE:
\be \label{eqinter2} Y^{n}_t =
\int_t^T f_n(s,Y^{n}_s) ds - \int_t^T Z^{n}_s dB_s, \,\,\,\, t \in [0,T].
\ee
The function $f_n$ being Lipschitz continuous w.r.t. $y$, uniformly in $t$,
therefore
through a result by Gobet et al. (2005) on numerical schemes of
BSDEs, this multidimensional equation  can be
numerically solved, at least in the case when the process $X$ is a
Markovian  diffusion. Therefore, this provides a way to simulate $Y^i$
since, as we will show it in Theorem 4 below, the sequence $(Y^{i,n})_{n\geq 0}$
converges to $Y^i$. Indeed, we have:
\bp For every $i \in {\cal J}$ and every $t \in [0,T]$,
the sequence $(Y^{i,n}_t)_{n \geq 0}$ is non-decreasing and P-a.s.
$Y^{i,n}_t \leq Y^i_t$.
\ep
\noindent{\it Proof}. For $n \in \N$,
and $k\in \N^*$, consider the following scheme. For every $i
\in {\cal J}$
$$
Y^{i,n,k}_t = \int_t^T \psi_i(s,X_s)ds +
n \int_t^T ( \max_{j\in {\cal J}^{-i}}(-\ell_{ij}(s)+Y^{j,n,k-1}_s)
- Y^{i,n,k}_s)^+ ds - \int_t^T Z^{i,n,k}_s dB_s,\,\,\,\, t \in [0,T]
$$
and
$$
Y^{i,n,0} =
E[\int_t^T\psi_i(s,X_s)ds|{\cal F}_t],\,\,\,t\leq T.
$$
From a result in El Karoui et al. (1997-1), $Y^{i,n,k}$
converges to $Y^{i,n}$ when $k$ tends to infinity. Now, let us show
by induction on $k$ that:
$$P-a.s.\,\,\,\,
Y^{i,n,k}_t \leq Y^{i,n+1,k}_t, \quad  n\geq 0,\,\,\,\,i \in {\cal J}, \,\,\,\,
t \in [0,T].
$$ For $k=0$ the property holds true.
Suppose now that it is also verified for some $k-1$ and let us show
that it is valid for $k$. For any $n\geq 0$, $i\in {\cal J}$ and
$t\in [0,T]$ we have:
$$
Y^{i,n+1,k}_t = \int_t^T \psi_i(s,X_s)ds +
(n+1) \int_t^T ( \max_{j\in {\cal
J}^{-i}}(-\ell_{ij}(s)+Y^{j,n+1,k-1}_s) - Y^{i,n+1,k}_s)^+ ds -
\int_t^T Z^{i,n+1,k}_s dB_s
$$
and
$$
Y^{i,n,k}_t = \int_t^T \psi_i(s,X_s)ds + n \int_t^T ( \max_{j\in {\cal
J}^{-i}} (-\ell_{ij}(s)+Y^{j,n,k-1}_s) - Y^{i,n,k}_s)^+ ds -
\int_t^T Z^{i,n,k}_s dB_s.
$$
Thanks to the induction hypothesis, for any $n\geq 0$, $i\in {\cal J}$ and
$t\leq T$, we have
$Y^{i,n,k-1}_t \leq Y^{i,n+1,k-1}_t$. Therefore, using the comparison
theorem of solutions of standard BSDEs (see e.g. El Karoui et al.
(1997-2), Theorem 2.2) we get that
$$
Y^{i,n+1,k}_t \geq Y^{i,n,k}_t,\quad t\leq T,
$$
which is the desired result. Now taking the limit as $k$ goes to
$+\infty$, we obtain, for any $n\geq 0$ and $i\in {\cal J}$, $Y^{i,n}
\leq Y^{i,n+1}$.

To finish the proof it remains to show that for any
$k\geq 0$ and $n\geq 0$ we have $Y^{i,n,k}_t \leq Y^i_t$ for any
$i\in {\cal J}$ and $t\leq T$ and then take the limit as $k$ goes
to infinity.
\noindent
Once more using induction on $k$, it hods true
that, for all $k\geq 0$, $n\geq 0$ and $ t\in [0,T]$,
$$
Y^{i,n,k}_t \leq Y^i_t,\quad
\mbox{ for any }i\in {\cal J}.
$$
Indeed, for $k=0$ the property is obviously satisfied. In order to go from $k$
to $k+1$, we
note that, by Eq.\@ (\ref{rBSDE1}), and since for that $k$, $Y^{i,n,k}_t \leq
Y^i_t$, for any $i\in {\cal J}$, it holds that for every $t\le T$, $(
\max_{j\in {\cal J}^{-i}}(-\ell_{ij}(t)+Y^{j,n,k}_t) - Y^{i}_t)^+=0.$

\noindent
Hence, for all $t\leq T$,
$$
Y^i_t= \int_t^T \psi_i(s,X_s)ds + n \int_t^T ( \max_{j\in {\cal J}^{-i}}
(-\ell_{ij}(s)+Y^{j,n,k}_s) - Y^{i}_s)^+ ds +K^i_T - K^i_t -
\int_t^T Z^{i}_s dB_s.
$$
Now, taking into account that the process $K^i$
is non-decreasing and finally and using the Comparison Theorem of
solutions of standard BSDEs, we get that
$$
Y^{i,n,k+1}_t \leq Y^i_t,\quad
\mbox{ for any }i\in {\cal J}.
$$

\noindent Finally taking the limit as
$k\rightarrow \infty$ we get that, for all $n\geq 0$ and $t\in [0,T]$,
$$
Y^{i,n}_t \leq Y^i_t, \quad\mbox{for any }i\in {\cal J}.
$$
The proof is now complete. $\Box$

\beth For any $i \in {\cal J}$ it holds true that:
$$E[\sup_{0\le t\leq T}|Y^{i,n}_t-Y^i_t|^2]\rightarrow 0 \mbox{ as
}n\rightarrow \infty.$$ \eeth

\noindent {\it Proof.} We have, for every $i \in {\cal J}$ and all $t \in
[0,T]$, $Y^{i,n}_t \leq Y^{i,n+1}_t$. Therefore there exists a
process $\bar{Y}^i$ such that,
$$
\lim_{n \to + \infty} Y^{i,n} = \bar{Y}^i_t \leq
Y^i_t,\quad t\in [0,T].
$$
Moreover, from (\ref{rBSDE2}) we get that, for any $t\leq
T$,
$$Y^{i,n}_t=\esssup_{\tau \geq t}E[\int_t^\tau \psi_i(s,X_s)ds + (L^{i,n}_\tau
\wedge Y^{i,n}_\tau)1_{[\tau <T]}|\cF_t].$$ This is due to the facts
that the process $n \int_0^. (L^{i,n}_s - Y^{i,n}_s)^+ ds$ is
increasing and satisfies $\int_0^T (Y^{i,n}_s - L^{i,n}_s \wedge
Y^{i,n}_s) n (L^{i,n}_s - Y^{i,n}_s)^+ ds=0$. Therefore, in order to
conclude, it is enough to use the representation result by El Karoui
et al. (1997-1) of solution of reflected BSDEs as Snell envelopes of
processes.

\noindent Now, since the process $Y^{i,n}_t + \int_0^t \psi_i(s,X_s)ds$ is a
continuous supermartingale, the non-decreasing limit
$\bar{Y}^i$ is a \cadlag process. Using now the result given in
Proposition \ref{thmsnell} - $v$, it follows that
$$
\bar{Y}^{i}_t=\esssup_{\tau \geq t}E[\int_t^\tau \psi_i(s,X_s)ds +
(\bar{L}^{i}_\tau \wedge \bar{Y}^{i}_\tau)1_{[\tau <T]}|\cF_t],\quad t\le T
$$
with $\displaystyle \bar{L}^{i}_t = \max_{k\in {\cal
J}^{-i}}(-\ell_{ik}(t)+\bar{Y}^{k}_t)$ is the nondecreasing limit of
$L^{i,n}$. But, from (\ref{rBSDE2}), taking expectation, dividing by
$n$ and taking the limit as $n\to \infty$ we obtain
$$
\int_0^T(\bar{L}^{i}_s-\bar{Y}^i_s)^+ds=0
$$
which implies that for any $t\leq T$, $\bar{Y}^i_t\geq \bar{L}^{i}_t$, since
these latter
processes are \cadlag. Now we can argue as in Section 4 to show that
the processes $\bar{Y}^i$, $i\in {\cal J}$, are continuous.
Therefore they satisfy the Verification Theorem whose solution is
unique. Hence, for any $i\in {\cal J}$, we have $\bar{Y}^i=Y^i$ and the
sequences $(Y^{i,n})_{n\geq 0}$ are nondecreasing and
converge to the continuous processes $Y^i$. Finally in order to
conclude we just need to use first Dini's Theorem and then the
Lebesgue dominated convergence theorem. $\Box$

\begin{rem}
It doesn't seem easy to obtain a convergence rate of $Y^{i,n}$ to
$Y^i$. In the two-modes case and when the switching costs are constant,
Hamad\`ene and Jeanblanc (2007)(Proposition 4.2) show that
the rate of convergence is $\frac{1}{n}$. This very interesting
issue will be addressed in a forthcoming work.
\end{rem}

\end{document}